\documentclass[11pt]{article}

\usepackage[english]{babel}

\usepackage[letterpaper,top=3cm,bottom=2.8cm,left=2.8cm,right=2.8cm,marginparwidth=1.75cm]{geometry}

\usepackage{amsmath}

\usepackage{amsthm}
\usepackage{amssymb}
\usepackage{mathrsfs}
\usepackage{tabularx}
\usepackage{amsfonts}
\usepackage{graphicx}
\usepackage{float}
\usepackage[colorlinks=true, allcolors=blue]{hyperref}
\usepackage{mathtools}
\usepackage{bbold}
\usepackage{indentfirst}
\usepackage{tikz}
\usetikzlibrary{intersections, calc, decorations.markings, arrows.meta}

\newcommand{\R}{\mathbb{R}}

\newcommand{\N}{{\mathbb{N}}}

\title{A new condition for the genericity of ergodic measures on Riemannian manifolds}

\author{Paul Mella}
\date{}

\newtheorem{theo}{Theorem}
\newtheorem{theoBis}{Theorem}

\newtheorem{prop}[theo]{Proposition}

\newtheorem{lemma}[theoBis]{Lemma}
\newtheorem{cor}[theo]{Corollary}

\theoremstyle{remark} 
\newtheorem*{rem}{Remark}

\begin{document}

\maketitle
\begin{abstract}
This article investigates the genericity of ergodic probability measures for the geodesic flow on Riemannian manifolds. We demonstrate that if the metric splits as a product metric within a tubular neighborhood of a geodesically complete submanifold containing a closed geodesic, then the closure of the set of ergodic measures does not encompass all invariant probability measures. Our findings notably provide an answer to the question of genericity of ergodic measures concerning a specific example of 3-manifold introduced by Gromov.
\end{abstract}

\section{Introduction}

\indent In the wake of a pioneering work of Hadamard dating back to 1898 \cite{H98}, the unit tangent bundle of a negatively curved Riemannian manifold $M$, which we will denote by $T^1M$ and equip with the geodesic flow, has been widely considered as a canonical example of hyperbolic dynamics. This paper focuses on the generalization of a property that has been known since the seventies under the negative curvature assumption: the genericity of ergodic measures in the space of invariant probability measures.

Let $\mathcal{M}$ be the space of invariant probability measures of the geodesic flow on $T^1M$, equipped with the weak-* topology induced by the geodesic distance of the Sasaki metric. It is a compact and convex subset of the space of probability measures, and its extremal points are exactly the ergodic measures, whose set will be denoted by $\mathcal M_e$. Hence, $\mathcal M_e$ is dense in $\mathcal M$ if and only if the latter is a realization of the Poulsen simplex \cite{LOS78}. Moreover, the set of ergodic probability measures $\mathcal{M}_e$ is a $G_\delta$ set of $\mathcal{M}$ \cite{P61}, thus the density condition $\overline{\mathcal{M}_e} = \mathcal{M}$ is equivalent to the following statement: \textit{ergodic measures are generic in the space of invariant probability measures}.

Since the works of Sigmund on topologically mixing Anosov flows \cite{S72}, this property has been known to hold for negatively curved compact manifolds. Coud\`ene and Schapira \cite{CS14} generalized it to non-positively curved manifolds whose universal cover has no flat strips. Their argument relies on several dynamical properties, such as the Closing Lemma, the transitivity and the local product structure. This bridge between the dynamical properties of a flow and the density of the set of ergodic measures has been generalized later by Gelfert and Kwietniak \cite{GK18}.

When it comes to necessary conditions, Coud\`ene and Schapira \cite{CS11, CS14} proved that for non-positively curved compact surfaces, the genericity of ergodic measures on $T^1M$ is equivalent to the condition that the universal cover of $M$ has no flat strips. For non-positively curved manifolds of higher dimension, Coud\`ene and Schapira remark that their argument could be generalized to prove that if $M$ contains an open isometric embedding of the product of $S^1$ with a Euclidean ball, then the genericity of ergodic measures fails. A similar setting was studied by Constantine, Lafont, McReynolds and Thompson \cite{CLMT19}. They introduced the notion of \textit{fat $1$-flat} as an open embedding of the product of $\R$ with a Euclidean ball. However, one of their main results is the existence of manifolds whose universal cover contains a \textit{fat $1$-flat}, yet the quotient of this flat is not isometric to the product of $S^1$ with a Euclidean ball. Proposition \ref{prop_gen} states a criterion that is sharper than the criterion of Coud\`ene and Schapira in dimension higher than 2. Notably, it applies to all manifolds whose universal cover contains a \textit{fat $1$-flat}, and does not require any general assumption on the curvature of $M$. 

A precise definition of tubular neighborhood can be found in \cite{L18}.

\begin{prop}\label{prop_gen}
Let $M$, $\widehat N$ and $\tilde N$ be Riemannian manifolds, $p:\tilde M\to M$ an isometric cover, and $\phi$ an open isometric embedding of $\widehat N\times \tilde N$, equipped with the product metric, into $\tilde M$. Assume that $\tilde N$ is geodesically complete and there exists $x\in \widehat N$ such that at least one of the following assumptions holds, with the notation $N=(p\circ\phi)(\{x\}\times \tilde N)$:
\begin{enumerate}
    \item\label{a1} $N$ is an embedded submanifold of $M$, $(p\circ\phi)(\widehat N\times \tilde N)$ is a tubular neighborhood of $N$ and both $N$ and $M\setminus \overline{N}$ contain a closed geodesic ($\overline N$ denotes the closure of $N$).
    \item\label{a2} $N$ is homeomorphic to $S^1$.
\end{enumerate}
Then, the set of ergodic probability measures is not generic in the set of invariant probability measures on $T^1M$.
\end{prop}

\begin{rem} The embedding $\phi$ being open implies $\dim \widehat N + \dim \tilde N=\dim M$.
\end{rem}


Proposition \ref{prop_gen} can be applied to the manifolds introduced by Constantine, Lafont, McReynolds and Thompson in \cite[Theorem A]{CLMT19}, with $N=S^1$ and $\widehat N$ a Euclidean ball, showing that ergodic measures are not dense in the space of invariant probability measures.

A direct corollary of Proposition \ref{prop_gen} is that the entropy density of $\mathcal M_e$ also fails under the same assumptions. Introduced by Orey and F{\"o}llmer \cite{FO88, O86}, this property has received a lot of interest in the last two decades, notably because of its relation to the weak specification property \cite{CLT20, QS16}. It has proven to be especially relevant in the field of symbolic dynamics \cite{PS07} and was used to prove large deviation theorems \cite{CTY17, EKW94}.

Section \ref{sec:proof} contains the proof of Proposition \ref{prop_gen}. In Section \ref{ex}, we state a corollary of Proposition \ref{prop_gen} and we present an example of a non-positively curved manifold that was introduced by Gromov \cite{G78} and has been studied by numerous authors \cite{BBE85, BGS85, E80, K98}. For this example, while the genericity of ergodic measures had remained an open question until now, we show that it fails using Corollary \ref{cor}.

\section{Proof of the non genericity of ergodic measures}
\label{sec:proof}

We consider $M, \tilde M, \widehat N, N, \tilde N$, $p$, $\phi$ and $x$ that satisfy the assumptions of Proposition \ref{prop_gen}. The goal of this section is to prove that $\mathcal M_e$ is not dense in $\mathcal M$. 

Without loss of generality, we can consider $\tilde N$ to be connected, since at least one of its connected components satisfies all the assumptions of Proposition \ref{prop_gen}. Moreover, if Assumption \ref{a1} does not hold, since $N$ is then compact, by making $\widehat N$ smaller if necessary, we can still assume that $(p\circ\phi)(\widehat N\times \tilde N)$ is a tubular neighborhood of $N$.  

Let us introduce some notations. Hereafter, $\pi:TM\rightarrow M$, $\tilde \pi:T\tilde M\rightarrow \tilde M$ and $\widehat{\pi}:T\widehat N\rightarrow \widehat N$ denote the projections that map any vector onto its base-point. The geodesic flows on the tangent bundles and unit tangent bundles of all the manifolds involved are denoted by the standard notation $(g_s)$, even when they are not complete. The Lebesgue measure on $\R$ is denoted by $\lambda$.

All the tangent bundles are equipped with the geodesic distance induced by the Sasaki metric. The distance to a subset of a manifold is defined as the infimum of all distances to the elements of this subset. It is a continuous map. For any $r>0$ and any subsets $A\subset T^1M$ and $\tilde A\subset T^1\tilde M$, their $r$-neighborhoods are denoted by: 
$$B_r(A) = \{v\in T^1M \mid d_{TM}(v, A) < r\}, \qquad \tilde B_r(\tilde A) = \{\tilde v\in T^1\tilde M \mid d_{T\tilde M}(\tilde v,\tilde A) < r\}$$

Similarly, the $r$-ball of $\widehat N$ centered at $x$ is denoted by $\widehat{B}_r = \{y\in \widehat N \mid d_{\widehat N}(x,y) < r\}$. According to \cite[Chapter 3, Proposition 4.2]{D92}, we can define a positive radius $\rho$ such that the ball $\widehat B_{\rho}$ is strongly convex, and making $\rho$ smaller if necessary, we can assume that the exponential map is well defined on the $\rho$-ball centered at $0$ of the tangent space $T_x\widehat N$. 

\medskip

Let $E$ be a connected component of $T^1N$. 

If $\dim N\ge 2$, then $E = T^1N$. According to Assumption \ref{a1}, there exists a closed geodesic in $M\setminus \overline N$. Hence, there exists a subset of $T^1M\setminus \overline E$ that is invariant by the geodesic flow and diffeomorphic to $S^1$: we denote it by $F$. Since $F$ is compact and disjoint from $\overline E$, the distance $d_{TM}(E, F) = \inf\{d_{TM}(x, y)\mid (x, y)\in E\times F\}$ is positive. 

If $\dim N\le 1$, then $N$ is homeomorphic to $S^1$, otherwise it would not contain a closed geodesic. Thus $T^1N$ has exactly two connected components, both compact. In this case, we set $F = T^1N\setminus E$, which still satisfies $d_{TM}(E, F)>0$. 

In both cases, we can pick $r\in (0, \frac\rho2)$ such that $2r < d_{TM}(E, F)$. The following Lemma describes the behavior of geodesics within the ball $\widehat B_{2r}$.

\begin{lemma}\label{lem1} Let $v\in T^1\widehat{B}_{r}\subset T^1\widehat N$. Then, $\{s>0 \mid \widehat{\pi}({g}_sv)\not \in \widehat{B}_{r}\}$ and $\{s>0 \mid \widehat{\pi}({g}_sv)\not \in \widehat{B}_{2r}\}$ are closed non-empty subsets of $\R$. Set:
$$s_1 = \min~\{s>0 \mid \widehat{\pi}({g}_sv)\not \in \widehat{B}_{r}\} \qquad \text{  and  } \qquad s_2 = \min~\{s>0 \mid \widehat{\pi}({g}_sv)\not \in \widehat{B}_{2r}\}$$
Then, $\{\widehat{\pi}({g}_sv) \mid s_1\le s< s_2\}$ is contained in $\widehat{B}_{2r}\setminus \widehat{B}_{r}$ and the inequality $s_1\leq \frac23s_2$ holds.
\end{lemma}

\begin{figure}[H]
\centering
\begin{tikzpicture}[scale=0.75]

    \coordinate (C) at (0,0);

    \draw[<->, >=stealth, dashed, very thick] 
        (C) -- ({1.6*cos(-45)}, {1*sin(-45)}) 
        node[midway, above right=-2pt] {$r$};

    \draw[<->, >=stealth, dashed, very thick] 
        (C) -- ({5*cos(-145)}, {3*sin(-145)}) 
        node[midway, below right=-2pt] {$2r$};

    \draw[very thick, name path=outer] (0,0) ellipse (5cm and 3cm);
    \draw[very thick, name path=inner] (0,0) ellipse (1.6cm and 1cm);

    \draw[very thick, name path=curve,
        decoration={markings, mark=at position 0.53 with {
            \draw[very thick, color=black!90!cyan, line join=round] (0.2, 0.15) -- (-0.05, 0) -- (0.2, -0.15);
        }},
        postaction={decorate}
    ] (-6.5, -0.2) .. controls (-3, 1) and (3, 1) .. (6.5, -0.4);

    \path [name intersections={of=curve and outer, sort by=curve, by={P3, P4}}]; 
    \path [name intersections={of=curve and inner, sort by=curve, by={P2, P_inner_right}}]; 

    \path[name path=vline] (0.8, -2) -- (0.8, 2);
    \path [name intersections={of=curve and vline, by={P1}}];

    \newcommand{\drawcross}[1]{
        \draw[very thick] (#1) ++(-0.18,0) -- ++(0.36,0);
        \draw[very thick] (#1) ++(0,-0.18) -- ++(0,0.36);
    }

    \drawcross{P3}
    \drawcross{P2}
    \drawcross{P1}
    \drawcross{C}

    \node[above left=2pt]  at (P3) {$\widehat{\pi}({g}_{s_2}v)$};
    \node[above left=2pt]  at (P2) {$\widehat{\pi}({g}_{s_1}v)$};
    \node[above=2pt]       at (P1) {$\widehat{\pi}(v)$};
    \node[left=4pt]        at (C)  {$x$};

\end{tikzpicture}
\caption{The behavior of geodesics within the ball $\widehat B_{2r}$.} 
\label{fig:temp}
\end{figure}

\begin{proof}Since the exponential map is well defined on the $\rho$-ball centered at $0$ of $T_x\widehat N$, the geodesic $s\mapsto \widehat \pi(g_sv)$ can be extended at least until it exits $\widehat B_{2r}$, thus $\{s>0 \mid \widehat{\pi}({g}_sv)\not \in \widehat{B}_{r}\}$ and $\{s>0 \mid \widehat{\pi}({g}_sv)\not \in \widehat{B}_{2r}\}$ are non-empty. They are closed because the balls $\widehat{B}_{r}$ and $\widehat{B}_{2r}$ are open.

Let $s \in [s_1,s_2)$. Then $\widehat{\pi}({g}_sv)$ belongs to $\widehat B_{2r}$ by definition of $s_2$. If it did belong to $\widehat{B}_{r}$, then, by convexity of $\widehat{B}_{r}$, there would exist a geodesic segment connecting $\widehat{\pi}(v)$ and $\widehat{\pi}({g}_sv)$ that is contained in $\widehat{B}_{r}$. By the strong convexity of $\widehat{B}_{\rho}$, this segment would have to be equal to $\{\widehat{\pi}({g}_tv) \mid 0\leq t\leq s\}$, which is contradictory because the latter is not contained in $\widehat{B}_{r}$. Hence, $\widehat{\pi}({g}_sv) \in \widehat{B}_{2r}\setminus \widehat{B}_{r}$. 

Moreover, by the strong convexity of $\widehat B_\rho$, the segment $\{\widehat \pi(g_tv)\mid 0\le t\le s_1\}$ is a length-minimizing geodesic, thus the triangle inequality yields $s_1\le 2r$. 

Finally, the reversed triangle inequality yields:
$$s_2-s_1=d_{\widehat N}\big(\widehat\pi(g_{s_2}v),\widehat\pi(g_{s_1}v)\big)\ge \big|d_{\widehat N}\big(x, \widehat\pi(g_{s_2}v)\big) - d_{\widehat N}\big(x, \widehat\pi(g_{s_1}v)\big)\big| \ge 2r-r = r$$

Hence, $s_1\le \frac23 s_2$.\end{proof}

\medskip

\begin{lemma}\label{lem2} Let $u\in T^1M\setminus B_{r}(E)$, and set $U = \left\{t>0 \mid g_tu \in B_r(E)\right\}$. Then,
$$\underset{T\rightarrow \infty}\liminf~\frac{ \lambda\big([0,T]\cap U\big)}{T} ~ \leq ~ \frac23$$
\end{lemma}

The proof of Lemma \ref{lem2} relies on the argument that the orbit of $u$ can be decomposed into a finite or countable union of geodesic arcs that intersect $E$ exactly once. The intersection with $E$ of those segments can be projected on geodesic segments of $\widehat N$, and Lemma \ref{lem1} yields an upper bound on the average time spent inside $B_r(E)$. It is worth noting that the estimate $\frac23$ in Lemma \ref{lem2} is arbitrary. Improving the argument of Lemma \ref{lem1} and choosing $r$ smaller could give arbitrarily low upper bounds.

\begin{proof}If $U$ is bounded, the limit of $\frac{ \lambda\left([0,T]\cap U\right)}{T}$ is zero and the Lemma is proved. Let us assume that $U$ is unbounded. As a non-empty open subset of $\R$, it is a finite or countable union of disjoint intervals.  Let $(t_{-1},{t_1})$ be one of those intervals (with $t_1 \in (0,\infty]$).

Let $\tilde E$ be a connected component of the unit tangent bundle of $\phi(\{x\}\times \tilde N)$, such that $E = dp(\tilde E)$. The geodesic distance induced by the Sasaki metric on $T\tilde M$ is bounded below by the geodesic distance on $\tilde M$. Hence, $\tilde\pi(\tilde B_r(\tilde E))$ is contained in the $r$-neighborhood of $\phi(\{x\}\times \tilde N)$, and since $\widehat N\times \tilde N$ is equipped with the product metric, the following inclusion holds:
\begin{equation}\label{lipschitz}\tag{$\ast$}\tilde \pi(\tilde B_r(\tilde E))\subset \phi\big(\widehat B_r\times \tilde N\big)\end{equation}

Since $\tilde N$ is connected and geodesically complete, $\phi(\{x\}\times \tilde N)$ is a connected component of $p^{-1}(N)$, and since $(p\circ \phi)(\widehat N\times \tilde N)$ is a tubular neighborhood of $N$, then $\phi(\widehat N\times \tilde N)$ is a connected component of $p^{-1}\big((p\circ\phi)(\widehat N\times \tilde N)\big)$. Thus any geodesic segment in $(p\circ\phi)(\widehat N\times \tilde N)$ admits a lift in $\phi(\widehat N\times \tilde N)$. In particular, we can find a lift $\tilde u\in T^1\tilde M$ of $u$ such that $g_{t}\tilde u$ belongs to $\tilde B_r(\tilde E)$ for all $t\in (t_{-1}, t_1)$. 

Let us pick $t_0 \in (t_{-1},{t_1})$, and let $v\in T\widehat N$ and $w\in T\tilde N$ such that $g_{t_0}\tilde u=d\phi(v,w)$. In order to prove by contradiction that $\|v\|$ is positive, let us assume that $\|v\|=0$. Then, $\tilde u = d\phi(v, g_{-t_0}w)$ holds, and one has $d_{T\widehat N}(v, 0_{T_x\widehat N}) = d_{\widehat N}(\widehat \pi(v), x)$. Hence,
\begin{equation}\label{zero_vector}\tag{$\ast\ast$}d_{T\tilde M}\big(g_{t_0}\tilde u, d\phi(0_{T_x\widehat N}, w)\big) ~=~ d_{T\tilde M}\big(\tilde u, d\phi(0_{T_x\widehat N}, g_{-t_0}w)\big) ~=~  d_{\widehat N}(\widehat \pi(v), x)\end{equation}

Since $g_{t_0}\tilde u\in\tilde B_r(\tilde E)$, Equation \eqref{lipschitz} shows that this distance is smaller than $r$ by \eqref{lipschitz}. Moreover, the vector $d\phi(0_{T_x\widehat N}, w)$ belongs to the unit tangent bundle of $\phi(\{x\}\times \tilde N)$. The latter is either equal to $\tilde E$ (if $\dim N\ge 2$) or the union of $\tilde E$ and a lift of $F$ (if $\dim N\le 1$). According to \eqref{zero_vector}, the vector $d\phi(0_{T_x\widehat N}, w)$ belongs also to the $2r$-neighborhood of $\tilde E$, which must be disjoint from any lift of $F$ since $2r<d_{TM}(E, F)$. In both cases, $d\phi(0_{T_x\widehat N}, w)$ must belong to $\tilde E$. Since $\tilde N$ is geodesically complete, the vector $d\phi(0_{T_x\widehat N}, g_{-t_0}w)$ also belongs to $\tilde E$, and \eqref{zero_vector} yields $\tilde u\in \tilde B_r(\tilde E)$, leading to a contradiction because $u\not\in B_r(E)$. Hence, the norm of $v$ must be positive. 

By \eqref{lipschitz}, since $g_{t_0}\tilde u\in \tilde B_r(\tilde E)$, the point $\widehat\pi\Big(\frac{v}{\|v\|}\Big)$ belongs to $\widehat B_r$, thus $\frac{v}{\|v\|}$ satisfies the assumptions of Lemma \ref{lem1} and we can set $s_1$ and $s_2$ as in the statement of the Lemma. Moreover, $g_s\frac{v}{\|v\|}$ is well defined for all $s\in [s_1, s_2)$, and the following holds:
$$\forall s \in [s_1, s_2), \qquad g_{\frac s{\|v\|}+t_0} \tilde u ~=~ {\|v\|}~d\phi\Big(g_s\frac v{\|v\|}, g_{s}\frac w{\|v\|}\Big)$$

Hence, the geodesic arc $\big\{\pi\big(g_{\frac s{\|v\|}+t_0}\tilde u\big)\mid s_1\le s<s_2\big\}$ is contained in $\phi(\widehat N\times\tilde N)$, which is a lift of a tubular neighborhood of $N$, and since the normal exponential map of $\phi(\{x\}\times \tilde N)$ coincides with he pull-back of the normal exponential map of $N$ by $dp$, one has:
$$\forall s \in [s_1, s_2), \qquad d_M\big(\pi\big(g_{\frac s{\|v\|}+t_0} u\big), N\big) ~=~  d_{\tilde M}\big(\tilde\pi\big(g_{\frac s{\|v\|}+t_0} \tilde u\big), \phi(\{x\}\times \tilde N)\big) ~\ge~ d_{\widehat N}\Big(\widehat \pi\Big(g_s\frac v{\|v\|}\Big), x\Big)$$

This is larger than $r$ by Lemma \ref{lem1}. Hence, the interval $\big[\frac{s_1}{\|v\|}+t_0, \frac{s_2}{\|v\|}+t_0\big)$ is disjoint from $U$, and the inequalities $\frac{s_1}{\|v\|}+t_0 \ge t_1$ and $\frac{s_2}{\|v\|}+t_0\le t_2$ hold, with $t_2=\inf\big(U\cap(t_1, +\infty)\big)$. Notice that $t_2$ exists because $U$ is unbounded and $t_1$ is not $+\infty$, since it is bounded by $\frac{s_1}{\|v\|}+t_0$. Additionally, the conclusion of Lemma \ref{lem1} yields:
$$t_1-t_0 ~\le~ \frac{s_1}{\|v\|} ~\le~ \frac{2s_2}{3\|v\|}~\le~ \frac23 (t_2 - t_0)$$

Since $t_1-t_{-1}=\sup\left\{t_1-t_0 \mid t_0\in(t_{-1},t_1)\right\}$, we also get the inequality:
$$t_1-t_{-1} ~\le~ \frac23 (t_2 - t_{-1})$$

\medskip

This argument can be applied to any maximal open interval in $U$. In particular, every such interval is bounded. Since $U$ is unbounded, its partition in open intervals must be countable, which means that there exist two sequences $\big(t_{-1}^{(n)}\big)_{n\in\N}$ and $\big(t_1^{(n)}\big)_{n\in\N}$ such that:

$$U = \underset{n\in\N}{\bigsqcup}\big(t_{-1}^{(n)},t_1^{(n)}\big)\quad\text{and}\qquad\forall n \in \N, \qquad\begin{cases}0<t_{-1}^{(n)}<t_1^{(n)}<t_{-1}^{(n+1)}\\
t_1^{(n)}-t_{-1}^{(n)}\le\frac23\big(t_{-1}^{(n+1)}-t_{-1}^{(n)}\big)
\end{cases}$$

Hence,
$$\forall n \ge1, \qquad \frac{\lambda\Big(\big[0,t_{-1}^{(n)}\big] \cap U\Big)}{t_{-1}^{(n)}} ~\leq~ \frac2{3t_{-1}^{(n)}}\sum_{k=0}^{n-1}\left(t_{-1}^{(k+1)} - t_{-1}^{(k)}\right)~\leq~\frac23$$

The sequence $t_{-1}^{(n)}$ tends to infinity because $U$ is unbounded, which proves Lemma \ref{lem2}.
\end{proof}

To conclude the proof of Proposition \ref{prop_gen}, let us build an invariant probability measure on $T^1M$ that is not in the closure of the set of ergodic probability measures.

\medskip

We can pick $\mu_E\in\mathcal M_e$ and $\mu_F\in\mathcal M_e$ respectively supported on $E$ and $F$ since they both contain a closed geodesic. Let us set $\mu = \frac34\mu_E + \frac14\mu_F\in\mathcal M$. 

The semi-norms induced by the bounded continuous functions form a basis of the weak-* topology. We define two bounded continuous functions $\chi_{r}:T^1M\to\R$ and $\xi_{r}:T^1M\to\R$ by:
$$\forall v\in T^1M\qquad \chi_{r}(v) = \frac{\max\big(0, r-d_{TM}(v, E)\big)}r\quad \text{and}\quad \xi_{r}(v) = \frac{\max\big(0, r-d_{TM}(v, F)\big)}r$$

The $r$-neighborhoods of $E$ and $F$ are disjoint since $2r<d_{TM}(E, F)$. Hence, the integrals of $\chi_r$ and $\xi_r$ are equal to $\frac {3}4$ and $\frac14$, respectively, and the following subset of $\mathcal M$ is an open neighborhood of $\mu$: 
$$\mathcal O = \left\{\nu\in\mathcal M\mid \int_{T^1M}\chi_{r}d\nu > \frac{2}3~~\text{and}~~\int_{T^1M}\xi_{r}d\nu > 0\right\}$$

By contradiction, assume that there exists $\nu\in\mathcal{O}\cap\mathcal M_e$. Since $\chi_{r}$ and $\xi_{r}$ are bounded by $1$ and vanish respectively on the complements of $B_r(E)$ and $B_r(F)$, the measure $\nu$ satisfies:
$$\nu(B_r(E))>\frac23\qquad\text{and}\qquad\nu(B_r(F))>0$$

According to the ergodic theorem, for $\nu$-almost all vectors $v\in T^1M$, one has:
\begin{equation*}\frac{\lambda\left(\left\{t\in[0,T] \mid g_tv \in B_r(E)\right\}\right)}T  ~~ \underset{T\rightarrow \infty}{\longrightarrow} ~~ \nu(B_r(E))>\frac23\end{equation*}

In particular, this limit holds for at least one vector $v\in B_r(F)$. This vector does not belong to $B_r(E)$ because of the inequality $2r<d_{TM}(E, F)$. Hence, Lemma \ref{lem2} can be applied to $v$, leading to a contradiction with the limit above. We have proved that the neighborhood $\mathcal O$ is disjoint from $\mathcal M_e$, hence the latter is not dense.

\section{Application to manifolds with a local product metric}

\label{ex}
In this section, we prove a direct corollary of Proposition \ref{prop_gen} and we use it to show that ergodic probability measures are not generic for a specific example of non-positively curved manifold.

\begin{cor}\label{cor}
Let $M$, $\widehat N$ and $N$ be Riemannian manifolds. Assume that $\dim \widehat N\ge 1$, $N$ is geodesically complete and contains a closed geodesic, and $M$ contains an open submanifold isometric to $\widehat N\times N$, equipped with the product metric. Then the set of ergodic probability measures is not generic in the set of invariant probability measures on $T^1M$.
\end{cor}

\begin{proof}Under those assumptions, there exists an open isometric embedding $\phi$ of $\widehat N\times \tilde N$ into $\tilde M$, where $\tilde M=M$, $\tilde N=N$, and the cover $p: \tilde M \to M$ is the identity map. Let $x\in\widehat N$ and let us prove that the submanifold $\phi(\{x\}\times N)$ satisfies Assumption \ref{a1} of Proposition \ref{prop_gen}.

Restricting, if necessary, $\widehat N$ to an open geodesic ball centered at $x$, we can assume that $\phi(\widehat N\times \tilde N)$ is a tubular neighborhood of $\phi(\{x\}\times N)$. Moreover, $\phi(\{x\}\times N)$ is geodesically complete and contains a closed geodesic since it is isometric to $N$. Finally, the set $\widehat N\setminus\{x\}$ is nonempty since $\dim \widehat N\ge 1$ and for all $y\in \widehat N\setminus\{x\}$, the submanifold $\phi(\{y\}\times N)\subset M\setminus \overline{\phi(\{x\}\times N)}$ also contains a closed geodesic. \end{proof}

Let us now discuss a specific example of a non-positively curved manifold that was introduced by Gromov \cite{G78} and has been studied by numerous authors  \cite{BBE85, BGS85, E80, K98}.

\medskip

Let $\Sigma_2$ be a compact hyperbolic surface of genus 2. Let $A\subset \Sigma_2$ be a simple closed geodesic such that $\Sigma_2\setminus A$ has two connected components isometric to each other. We denote by $\Sigma_{\mathrm{half}}$ one of those connected components. Rescaling $\Sigma_2$ if necessary, let us assume that $A$ is isometric to $S^1$. Let $\phi$ be that isometry. Then, $(\Sigma_{\mathrm{half}} \cup A)\times S^1$ and $S^1 \times (\Sigma_{\mathrm{half}} \cup A)$ are compact 3-manifolds with boundary, and the following is an isometry between their boundary components:
$$\begin{aligned}
\Phi \colon \partial\left((\Sigma_{\mathrm{half}} \cup A)\times S^1\right) &\to \partial\left(S^1 \times (\Sigma_{\mathrm{half}} \cup A)\right) \\
(x, \theta) &\mapsto (\phi(x), \phi^{-1}(x))
\end{aligned}$$

\begin{figure}[H]
\centering
\includegraphics[angle=0, scale=0.35]{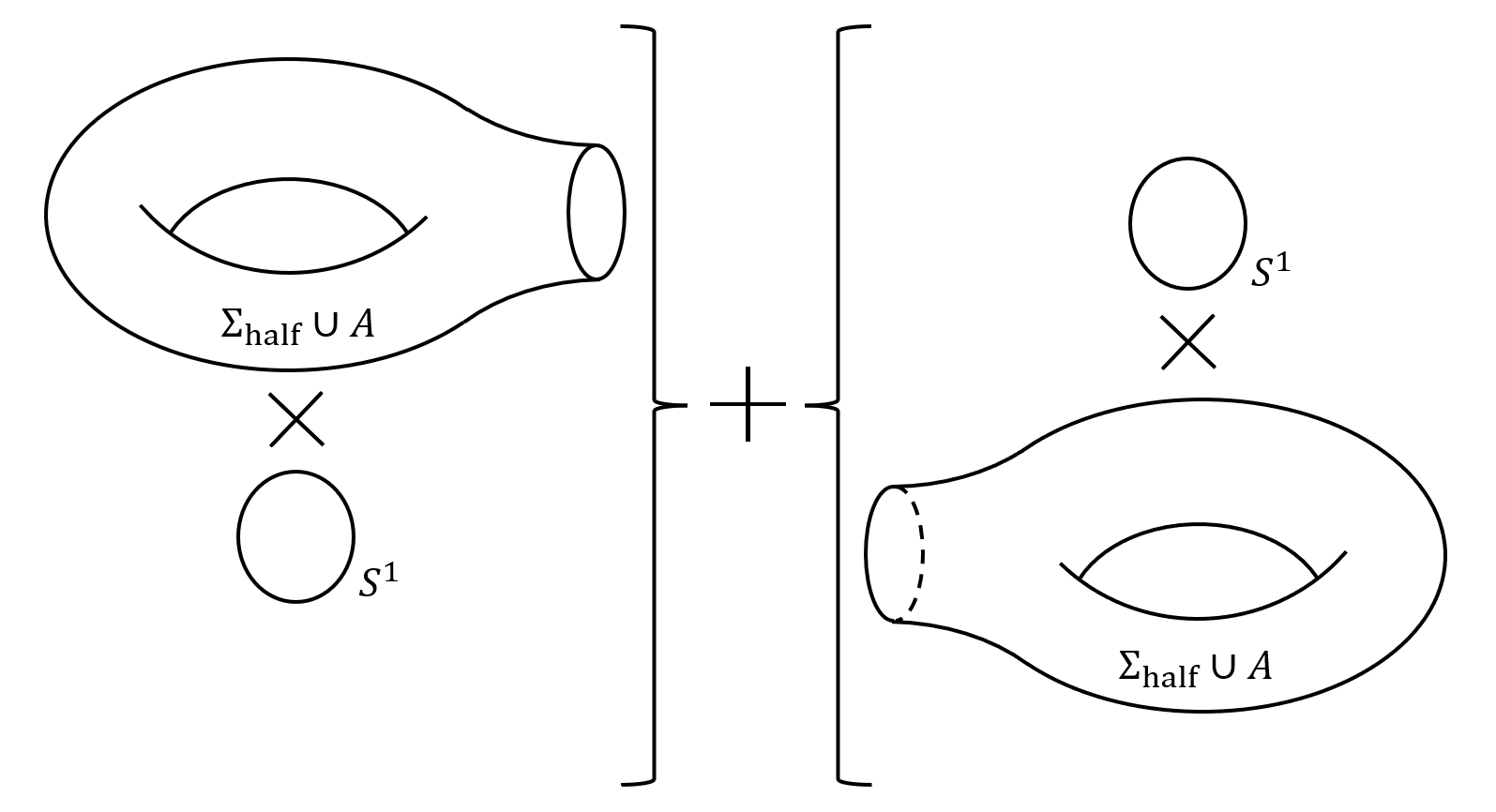}
\caption{The non-positively curved $3$-manifold introduced by Gromov.} 
\label{fig:Manifold_Gromov}
\end{figure}

Finally, modify the constant-curvature Riemannian metric of $\Sigma_{\mathrm{half}}$ so that the adjunction space
$$M = \left((\Sigma_{\mathrm{half}} \cup A)\times S^1\right) \sqcup_{\Phi} \left(S^1\times (\Sigma_{\mathrm{half}} \cup A)\right)$$
is a smooth $3$-manifold whose curvature vanishes on the central torus. This construction is illustrated in Figure \ref{fig:Manifold_Gromov}. 

Even though the universal cover of $M$ does not contain a fat $1$-flat, we can prove that $\mathcal M_e$ is not dense in $\mathcal M$.

\begin{cor} Ergodic measures are not generic in the space of invariant probability measures on $T^1M$.
\end{cor}

\begin{proof}
The inclusion map $\Sigma_{\mathrm{half}}\times S^1 \rightarrow M$ is an embedding that satisfies the assumptions of Corollary \ref{cor}, with $S^1$ being a geodesically complete manifold that contains at least one closed geodesic.
\end{proof}

\medskip

\bibliographystyle{plain} 
\bibliography{sample.bib}

\end{document}